\documentclass{article}

\usepackage{amssymb}

\usepackage[english,francais]{babel}

\newtheorem{e-proposition}[theorem]{Proposition}

\newtheorem{e-definition}[theorem]{Definition\rm}

\newtheorem{theoreme}{Th\'eor\`eme}
\newtheorem{lemme}[theoreme]{Lemme}

\newtheorem{definition}[theoreme]{D\'efinition\rm}

\def\og{\leavevmode\raise.3ex\hbox{$\scriptscriptstyle\langle\!\langle$~}}
\def\fg{\leavevmode\raise.3ex\hbox{~$\!\scriptscriptstyle\,\rangle\!\rangle$}}

\usepackage{amsmath}
\usepackage{a4wide}
\newcommand{\cqfd}{~\hspace{\stretch{1}}$\square$}

\date{\today}

\begin{document}
%
\selectlanguage{francais}
\title{Remarques \`a propos de l'op\'erateur de Dirac cubique}
\author{Nicolas Prudhon}
\maketitle
\medskip
\selectlanguage{francais}

\begin{abstract}
\selectlanguage{francais}
		 En 1999, Kostant 
introduit un op\'erateur de Dirac cubique $D_{\mathfrak{g}/\mathfrak{h}}$
associ\'e \`a tout triplet $(\mathfrak{g},\mathfrak{h},B)$, o\`u
$\mathfrak{g}$ est une alg\`ebre de Lie complexe munie de la forme bilin\'eaire sym\'etrique $\mathop{ad}\mathfrak{g}$-invariante non d\'eg\'en\'er\'ee
$B$, et $\mathfrak{h}$ est une sous-alg\`ebre de Lie de $\mathfrak{g}$ sur laquelle $B$ est non d\'eg\'en\'er\'ee.
Kostant montre alors que le carr\'e de $D_{\mathfrak{g}/\mathfrak{h}}$ 
v\'erifie une formule qui g\'en\'eralise la formule de Parthasarathy. 
Nous donnons ici une nouvelle d\'emonstration de cette formule. Tout d'abord, au moyen d'une
induction par \'etage, nous montrons qu'il suffit d'\'etablir la formule dans le cas particulier 
o\`u $\mathfrak{h}={0}$. Il appara\^it alors que, dans ce cas, l'annulation du terme d'ordre $1$ dans la formule de 
Kostant pour $D_{\mathfrak{g}/\mathfrak{h}}^2$ est une cons\'equence de propri\'et\'es classiques en cohomologie des alg\`ebres de Lie, tandis que le fait que le carr\'e du terme cubique soit scalaire r\'esulte de telles consid\'erations, 
ainsi que de l'identit\'e de Jacobi.

{\it Pour citer cet article~: N. Prudhon, C. R. Acad. Sci. Paris, Ser. I *** (20**).}
\vskip 0.5\baselineskip

\selectlanguage{english}
\noindent{\bf Abstract}
\vskip 0.5\baselineskip
\noindent
{\bf Remarks on the Kostant Dirac operator}
In 1999, Kostant \cite{Kostant1999} indroduces a Dirac operator $D_{\mathfrak{g}/\mathfrak{h}}$
associated to any triple $(\mathfrak{g},\mathfrak{h},B)$, where $\mathfrak{g}$ is a complex Lie algebra
provided with an $\mathop{ad}\mathfrak{g}$-invariant non degenerate nsymetric bilinear form $B$, 
and $\mathfrak{h}$ is a Lie 
subalgebra of $\mathfrak{g}$ such that the bilinear form $B$ is non degenerate on $\mathfrak{h}$. 
Kostant then shows that the square of this operator safisties a formula that generalizes 
the so-called Parthasarathy formula \cite{Parthasarathy1972}. We give here a new proof of this formula.
First we use an induction by stage argument to reduce the proof of the formula to the particular case
where $\mathfrak{h}=0$. In this case we show that the vanishing of the first ordrer term
in the Kostant formula for $D_{\mathfrak{g}/\mathfrak{h}}^2$ is a consequence
of classic properties related to Lie algebra cohomology, and the fact that the square of the cubic term
is a scalar follows from such considerations, together with the Jacobi identity.
{\it To cite this article: N. Prudhon, C. R. Acad. Sci. Paris, Ser. I *** (20**).}
\end{abstract}


\bigskip
\selectlanguage{francais}

Commen\c cons par d\'efinir l'op\'erateur de Dirac cubique.
Soit $(\mathfrak{g},\mathfrak{h},B)$ un triplet o\`u $\mathfrak{g}$ est une alg\`ebre de Lie 
complexe, $\mathfrak{h}$ une sous-alg\`ebre de Lie de $\mathfrak{g}$, et $B$ est une forme
bilineaire sym\'etrique sur $\mathfrak{g}$ invariante par l'action adjointe, i.e.
\begin{equation}\label{ad-invariance}
 B([x,y],z)+B(y,[x,z])=0\,,\quad (\forall x,y,z \in \mathfrak{g})\,.
\end{equation}
Nous supposons en outre que $B$ est non d\'eg\'en\'er\'ee, ainsi que sa restriction \`a $\mathfrak{h}$.
La restriction de $B$ \`a l'orthogonal $\mathfrak{h}^\perp$ de $\mathfrak{h}$ pour $B$, est alors
non d\'eg\'en\'er\'ee, et on obtient la d\'ecomposition $\mathop{ad}\mathfrak{l}$-invariante
$\mathfrak{g}=\mathfrak{h}\oplus\mathfrak{h}^\perp\,. $
Ainsi l'alg\`ebre $\mathfrak{h}$ est repr\'esent\'ee dans $\mathfrak{h}^\perp$ par restriction de l'action adjointe, 
induisant un homomorphisme
$$	\nu \colon \mathfrak{h} \to \mathfrak{so}(\mathfrak{h}^\perp)\,.$$
Comme la forme bilin\'eaire $B$ est non d\'eg\'en\'er\'ee sur $\mathfrak{h}^\perp$, nous pouvons \'egalement d\'efinir
l'alg\`ebre de Clifford $\mathcal{C}(\mathfrak{h}^\perp)$. Plus pr\'ecis\'ement, l'alg\`ebre de Clifford
est le quotient l'alg\`ebre tensorielle $T(\mathfrak{h}^\perp)$ par l'id\'eal engendr\'e par les \'elements de la forme
\begin{equation}\label{relation-clifford}
	xy+yx-2B(x,y)\,,\text{ o\`u } x,y \in \mathfrak{h}^\perp\,.
\end{equation}
Nous identifierons l'espace vectoriel sous-jacent \`a l'alg\`ebre de Clifford $\mathcal{C}(\mathfrak{h}^\perp)$ \`a celui de 
l'alg\`ebre ext\'erieur $\wedge \mathfrak{h^\perp}$
au moyen de l'isomorphisme de Chevalley. 
Pour tout $u,v \in\wedge \mathfrak{h^\perp}$,
nous distinguerons alors le produit de Clifford $uv$ et le produit ext\'erieur $u\wedge v$.
Le produit ext\'erieur par $u$ sera not\'e $e(u)$.
La forme $B$ s'\'etend \`a $\wedge \mathfrak{h^\perp}$ en une forme
non d\'eg\'en\'er\'ee et fournit ainsi un
isomorphisme de l'alg\`ebre ext\'erieure avec son dual : $\wedge \mathfrak{h^\perp}=\left( \wedge \mathfrak{h^\perp}\right)^*$.
La transpos\'ee $\iota(x)$ de $e(x)$ pour $x\in\mathfrak{h^\perp}$, peut alors \^etre vue comme endomorphisme de 
$\wedge \mathfrak{h^\perp}$, et l'on a alors, pour $w\in \wedge\mathfrak{h^\perp}$ et $x\mathfrak{h}^\perp$,
$$	xw=\big(e(x)+\iota(x)\big)w\,.$$
Les relations de Clifford (\ref{relation-clifford}) s'\'etendent \'egalement de la fa\c con suivante :
\begin{equation}\label{relation-clifford-2}
	xw-(-1)^kwx = 2\iota(x)w \,,\text{ o\`u } x\in \mathfrak{h}^\perp \text{ et } w\in\wedge^k \mathfrak{h}^\perp\,.
\end{equation}
En outre, \`a travers ces isomorphismes, la restriction de la $3$-forme fondamentale de $\mathfrak{g}$ \`a $\mathfrak{h}^\perp$, i.e.
$$	v \colon \wedge^3 \mathfrak{h}^\perp \to \mathbb{C}\,,\quad B(v,x\wedge y\wedge z)=-\frac{1}{2}B(x,[y,z])\,,$$
d\'efinit un \'el\'ement $v_{\mathfrak{g}/\mathfrak{h}}\in \mathcal{C}(\mathfrak{h}^\perp)$.
Soit maintenant $(X_i)$ une base orthonorm\'ee de $\mathfrak{h}^\perp$, et $\mathcal{U}(\mathfrak{g})$ l'alg\`ebre de Lie
enveloppante de $\mathfrak{g}$.
\begin{definition}
L'op\'erateur de Dirac cubique est l'\'el\'ement
$D_{\mathfrak{g}/\mathfrak{h}} \in \mathcal{U}(\mathfrak{g})\otimes\mathcal{C}(\mathfrak{h}^\perp)$ d\'efini par
$$	D_{\mathfrak{g}/\mathfrak{h}} = \sum_i X_i\otimes X_i + 1\otimes v_{\mathfrak{g}/\mathfrak{h}}\,.$$
\end{definition}
Le qualificatif cubique exprime le fait que $v_{\mathfrak{g}/\mathfrak{h}}$ est de degr\'e $3$.

Il existe un unique homomorphisme d'alg\`ebre $\Delta_\mathfrak{h}$ d\'efini par
$$	\Delta_\mathfrak{h} \colon \mathcal{U}(\mathfrak{h}) \to \mathcal{U}(\mathfrak{g})\otimes\mathcal{C}(\mathfrak{h}^\perp)\,,
	\qquad \Delta_\mathfrak{h}(y)=y\otimes 1 + 1 \otimes \nu(y) \quad (y\in \mathfrak{h})\,.$$
La $\mathbb{Z}$-graduation sur l'alg\`ebre tensorielle induit une $\mathbb{Z}_2$-graduation sur l'alg\`ebre de Clifford. En effet,
l'automorphisme $\kappa$ de $\mathfrak{h}^\perp$ d\'efini par $\kappa(x)=-x$ v\'erifie $\kappa(x)^2=1$, et se prolonge 
donc d'apr\`es (\ref{relation-clifford})
en un automorphisme de $\mathcal{C}(\mathfrak{h}^\perp)$ encore not\'e $\kappa$. La $\mathbb{Z}_2$-graduation est alors
donn\'ee par la d\'ecomposition de $\mathcal{C}(\mathfrak{h}^\perp)$ en sous-espaces propres. 
En notant $\mathop{\otimes}\limits^{\_}$ le produit tensoriel gradu\'e, 
nous avons un isomorphisme d'alg\`ebres gradu\'ees
$$	\mathcal{C}(\mathfrak{g}) = \mathcal{C}(\mathfrak{h}^\perp) \mathop{\otimes}\limits^{\_} \mathcal{C}(\mathfrak{h})\,.$$
En consid\'erant la graduation triviale sur $\mathcal{U}(\mathfrak{g})$, on obtient une $\mathbb{Z}_2$-graduation
sur le produit tensoriel $\mathcal{U}(\mathfrak{g}) \otimes \mathcal{C}(\mathfrak{h}^\perp)$.
\begin{lemme}\label{h-invariance}
L'op\'erateur de Dirac cubique est alors $\mathfrak{h}_\Delta$-invariant, i.e.
$$	
D_{\mathfrak{g}/\mathfrak{h}} \in 
\big( \mathcal{U}(\mathfrak{g}) \otimes \mathcal{C}(\mathfrak{h}^\perp) \big)^{\mathfrak{h}_\Delta}\,,
$$
commute, au sens gradu\'e, avec l'image $\mathfrak{h}_\Delta=\Delta_{\mathfrak{h}}(\mathfrak{h})$.
\end{lemme}
Soit $(Y_j)$ une base orthonorm\'ee de $\mathfrak{h}$, $\Omega_\mathfrak{h}=\sum_jY_j^2 \in \mathcal{U}(\mathfrak{h})$ l'op\'erateur
de Casimir pour $\mathfrak{h}$, et $\Omega_\mathfrak{g}=\sum_iX_i^2 + \sum_jY_j^2 \in \mathcal{U}(\mathfrak{g})$ l'op\'erateur
de Casimir pour $\mathfrak{g}$.
\begin{theoreme}\cite{Kostant1999} \label{Kostant}
Le carr\'e de $v_{\mathfrak{g}/\mathfrak{h}}$ dans l'alg\`ebre de Clifford est 
un scalaire $c_{\mathfrak{g}/\mathfrak{h}}$. De plus,
$$	
D_{\mathfrak{g}/\mathfrak{h}}^2 = 
\Omega_\mathfrak{g}\otimes 1 - \Delta_\mathfrak{h}(\Omega_\mathfrak{h})+c_{\mathfrak{g}/\mathfrak{h}}\, 
$$
\end{theoreme}

Le premier argument de la d\'emonstration que nous donnons ici, de type induction par \'etage,
appara\^it dans \cite{HuangPandzic2006,HuangPandzic-book}. On pourra \'egalement consulter
\cite{MehdiZierau2006}. Notons $D_\mathfrak{r}$ l'op\'erateur de Dirac Kostant associ\'e \`a un triplet 
de la forme $(\mathfrak{r},0,B)$.
\begin{lemme} \cite[Theorem 9.4.1]{HuangPandzic-book} \\
 (i) On a la d\'ecomposition suivante
$$	
D_\mathfrak{g} = D_{\mathfrak{g}/\mathfrak{h}}\mathop{\otimes}\limits^{\_}1 + 
\Delta_\mathfrak{h}\mathop{\otimes}\limits^{\_}1(D_\mathfrak{h}) \,.
$$
 (ii) Les composantes $D_{\mathfrak{g}/\mathfrak{h}}\mathop{\otimes}\limits^{\_}1$ et 
$\Delta_\mathfrak{h}\mathop{\otimes}\limits^{\_}1(D_\mathfrak{h})$ anticommutent. 
\end{lemme}
Ce lemme r\'esulte essentiellement du lemme \ref{h-invariance} et de la d\'efinition du produit tensoriel gradu\'e.

La cons\'equence de ce lemme dont nous avons besoin est la suivante :
$$	
D_{\mathfrak{g}/\mathfrak{h}}^2\mathop{\otimes}\limits^{\_}1=
D_\mathfrak{g}^2-\Delta_\mathfrak{h}\mathop{\otimes}\limits^{\_}1\left(D_{\mathfrak{h}}^2\right)\,.
$$
Ainsi, le th\'eor\`eme \ref{Kostant} r\'esulte du m\^eme th\'eor\`eme pour les alg\`ebres de Lie quadratiques,
i.e. les triplets du type $(\mathfrak{r},0,B)$. En outre,
on a \'egalement 
$$	c_{\mathfrak{g}/\mathfrak{h}}=c_{\mathfrak{g}}-c_{\mathfrak{h}}\,.$$
Il reste maintenant \`a montrer le th\'eor\`eme dans le cas o\`u $\mathfrak{h}=0$
(et donc $\mathfrak{h}^\perp=\mathfrak{g}$).
\begin{theoreme}
Soit $(\mathfrak{g},B)$ une alg\`ebre de Lie quadratique.
Le carr\'e dans l'alg\`ebre Clifford $\mathcal{C}(\mathfrak{g})$ de la $3$-forme fondamentale 
est un scalaire $c_\mathfrak{g}$ et 
$$	D_\mathfrak{g}^2 = \Omega_\mathfrak{g} + c_\mathfrak{g}\,.$$
\end{theoreme}

\noindent \emph{D\'emonstration.}
Posons $v=v_\mathfrak{g}$ pour simplifier les notations et calculons $D_\mathfrak{g}^2$. Pour ceci,
introduisons la transpos\'ee du corchet de Lie $\delta \colon \mathfrak{g} \to \wedge^2\mathfrak{g}$. Nous avons
$$	\delta(x)=\sum_{i<j}B(x,[X_i,X_j])X_i\wedge X_j\,,\quad(x\in\mathfrak{g})\,.$$
Par la suite, $\wedge^2\mathfrak{g}$ est identifi\'e \`a un sous-espace de $\mathcal{C}(\mathfrak{g})$ comme
pr\'ec\'edement, de sorte que $\delta(X) \in \mathcal{C}(\mathfrak{g})$.
Introduisons \'egalement la d\'erivation $d_v$ de l'alg\`ebre gradu\'ee $\mathcal{C}(\mathfrak{g})$ donn\'ee par 
$$ d_v(a)=va-\kappa(a)v\,,\quad(a\in\mathcal{C}(\mathfrak{g})\,.$$
Remarquons que $d_v$ est bien une d\'erivation car
$$	d_v(ab)=vab-\kappa(ab)v=vab-\kappa(a)\kappa(b)v=(va-\kappa(a)v)b+\kappa(a)(vb-\kappa(b)v)=(d_va)b+\kappa(a)d_vb\,.$$
Il vient,
\begin{align*}
D_\mathfrak{g}^2&= \left(\sum_i X_i\otimes X_i\right)\left( \sum_j X_j\otimes X_j\right)
		   + \sum_k X_k\otimes d_v(X_k) + v^2 \\
		&= \left(\sum X_i^2\right)\otimes1 + \sum_{i<j}[X_i,X_j]\otimes X_iX_j
		   + \sum_k X_k\otimes d_v(X_k) + v^2
\end{align*}
Or,
$$ 
\sum_{i<j}[X_i,X_j]\otimes X_iX_j = \sum_k X_k \otimes \left( \sum_{i<j} B(X_k,[X_i,X_j])X_iX_j \right)
=  \sum X_k\otimes \delta(X_k)\,.$$
Ainsi,
$$	D_\mathfrak{g}^2=\Omega_\mathfrak{g}\otimes 1+v^2+ \sum_k X_k\otimes (\delta+d_v)(X_k)\,.$$ 	
Il suffit donc de montrer que $\delta+d_v=0$ et que $v^2$ est scalaire. 
La relation $\delta+d_v=0$ peut \^etre \emph{v\'erifi\'ee} en quelques lignes
en \'ecrivant $d_v$ explicitement dans une base orthonorm\'ee \cite[Lemma 3.3]{Agricola2003}. Nous pr\'ef\`erons la d\'emonstration
suivante, qui explique davantage \emph{pourquoi} une telle relation est vraie.
 
Pour cela commen\c cons par \'etendre $\delta$ en une d\'erivation de l'alg\`ebre $\wedge\mathfrak{g}$.
Pour $X\in \mathfrak{g}$, notons $X^*=B(X,\cdot)$ la forme lin\'eaire associ\'ee. On peut alors
\'ecrire $\delta(X^*)(Y,Z)=X^*([Y,Z])$. Soit $C^k(\mathfrak{g},\mathbb{C})$ l'espace des applications $k$-lin\'eaires
sur $\mathfrak{g}$, non n\'ecessairement altern\'ees. Alors $\delta$ n'est autre que la restriction de la diff\'erentielle
$$	d \colon C^k(\mathfrak{g},\mathbb{C}) \to C^{k+1}(\mathfrak{g},\mathbb{C}) $$ 
donn\'ee par
$$ 
dw(x_0,\ldots,x_k)=
\sum_{s<t} (-1)^{s}w(x_0,\ldots,\hat{x}_s,\ldots,[x_s,x_t],\ldots,x_k)\,.
$$
Le lien fondamental entre $v$, $\delta$ et $B$ est alors donn\'e par la relation $dB=2v$. En effet, comme 
la forme bilin\'eaire $B$ est $\mathop{ad}{\mathfrak{g}}$-invariante, 
$$ dB(x,y,z)=B([x,y],z)+B(y,[x,z])-B(x,[y,z])=-B(x,[y,z])=2v(x,y,z)\,.$$
Pour $X\in\mathfrak{g}$, notons $\theta_X$ l'action de $X$ sur $C^k(\mathfrak{g},\mathbb{C})$,
$$	\theta_Xw(x_1,\ldots,x_k)=\sum_s w(x_1,\ldots,[X,x_s],\ldots,x_k)\,,$$
Remarquons que l'invariance de $B$, \'equation~(\ref{ad-invariance}), s'\'ecrit alors $\theta_XB=0$.
Rappelons encore que $\iota_X$ est d\'efini sur $C^\bullet(\mathfrak{g},\mathbb{C})$ par $\iota_Xw(\cdot)=w(X,\cdot)$. 

On a \'egalement la formule de Cartan
$$	\iota_Xd+d\iota_X=\theta_X\,.$$
Ainsi, pour $X \in \mathfrak{g}$, utilisant l'\'equation~(\ref{relation-clifford-2}), 
$$	d_v(X)=2\iota_Xv=\iota_XdB=-d\iota_XB+\theta_XB=-dX^*=-\delta(X)\,.$$
Nous avons donc obtenu la relation recherch\'ee $\delta+d_v=0$.

Venons-en maintenant \`a l'examen de $v^2$.
L'op\'erateur $d$
laisse stable l'espace des formes altern\'ees $\wedge \mathfrak{g}^*$, alg\`ebre sur lequel il d\'efinit une d\'erivation
- pour la graduation d\'efinie par la parit\'e du degr\'e. 
Nous venons donc de montrer que $d$ et $d_v$ \'etaient des d\'erivations, qui \`a une consante pr\`es, co\"incidaient
en degr\'e $1$.
Il en est donc de m\^eme sur toute l'alg\`ebre $\wedge \mathfrak{g}$. 
Le fait que $d^2=0$ est une cons\'equence imm\'ediate de l'identit\'e de Jacobi. En effet, en degr\'e $1$,
$$ d^2w(x,y,z)=dw([x,y],z)+dw(y,[x,z])-dw(x,[y,z])=w([[x,y],z]+[y,[x,z]]-[x,[y,z]])=0\,.$$
Autrement dit, toujours en degr\'e $1$, $d^2=(\delta\otimes1)\circ \mathop{Alt} \circ\, \delta = 0$,
car le coproduit $\delta$ v\'erifie l'identit\'e de Jacobi duale. 
Nous en d\'eduisons que $d_v^2=0$. Par ailleurs, comme l'automorphisme $\kappa$ v\'erifie $\kappa(v)=-v$ et $\kappa^2(a)=a$,
$$	d_v^2a=v(va-\kappa(a)v)-\kappa(va-\kappa(a)v)v=v^2a-av^2\,.$$
Par cons\'equent, $v^2$ est dans le centre de l'alg\`ebre de Clifford. En dimension paire, le centre
est r\'eduit au scalaire, et donc $v^2$ est scalaire. Dans le cas de la dimension impaire, le centre
est engendr\'e par les scalaires et les \'el\'ements de degr\'e maximal, qui sont donc de degr\'e impair. 
Or, comme $v$ est de degr\'e $3$, son carr\'e $v^2$ ne
peut contenir que des termes de degr\'e pair. Donc $v^2$ est scalaire dans ce cas \'egalement.\cqfd


\end{document}